\documentclass{svmult}
\usepackage{mathrsfs}
\usepackage{amssymb}
\newenvironment{enumerate-roman}{\begin{enumerate}}{\end{enumerate}}
\usepackage{makeidx}         
\usepackage{graphicx}        
\usepackage{multicol}        
\usepackage[bottom]{footmisc}

\makeindex             

\def\proof{{\bf Proof:} }
\def\endpf{\begin{flushright} \hfill $\Box$
                   \end{flushright} }
\def\p{\frac{\partial}{\partial x}}
\def\w{W_{\nu\Delta}}
\def\wk{W_{\nu\Delta,k}}

\begin{document}

\title*{Pathwise stationary solutions of stochastic Burgers equations with $L^2[0,1]$-noise
and stochastic Burgers integral equations on infinite horizon}
\titlerunning{Stationary Solution of stochastic Burgers equation}
\author{Yong Liu\inst{1,2},
Huaizhong Zhao\inst{1}}
\authorrunning{Y. Liu and H.Z. Zhao}
\institute{ Department of Mathematical Sciences, Loughborough
University, Loughborough, LE11 3TU, UK.
\and
LMAM, School of Mathematical Sciences,
                      Peking University, Beijing 100871,P.R. China
                      \\
\texttt{liuyong@math.pku.edu.cn}, \texttt{H.Zhao@lboro.ac.uk} }
\date{}
\maketitle
\begin{abstract}
In this paper, we show the existence and uniqueness of the
stationary solution $u(t,\omega)$ and stationary point $Y(\omega)$ of
the differentiable random dynamical system $U:R\times L^2[0,1]\times \Omega\to
L^2[0,1]$ generated by  the
stochastic Burgers equation with $L^2[0,1]$-noise and large
viscosity, especially, $u(t,\omega)=U(t,Y(\omega),\omega)=Y(\theta(t,\omega))$, and $Y(\omega)
\in H^1[0,1]$
 is the unique solution of the following equation in $L^2[0,1]$
$$
  Y(\omega)=\frac{1}{2}\int_{-\infty}^0T_\nu(-s)\frac{\partial
  (Y(\theta(s,\omega))^2}{\partial
  x}ds+\int_{-\infty}^0T_\nu(-s)dW_s(\omega),
$$
where $\theta$ is the group of $P$-preserving ergodic transformation
on the canonical probability space $(\Omega, {\cal F}, P)$ such that
$\theta(t,\omega)(s)=W(t+s)-W(t)$.
\end{abstract}
\vskip-10pt


\vskip0.5cm
\section{Introduction}
The stationary point (or stationary solution) is one of the
fundamental concepts in dynamical systems. For example, for an
autonomous ordinary differential equation (ODE): $\dot{X}=F(X)$, a
stationary point is a point in the set $\{x: F(x)=0\}$ in the phase
space or the stationary solution is a trajectory (fixed in the
autonomous case) satisfying $F(X(t))=0$ for any $t\in
(-\infty,\infty)$. Roughly speaking, the behaviour of the solution
near a stationary point describes the asymptotic properties of the
dynamical systems and the stationary solution gives
equilibrium state. For the infinite-dimensional dynamical systems
generated by some partial differential equations (PDEs) of the
following form (see \cite{[R]}, \cite{[T1]})
$$
 \frac{\partial u(t,x)}{\partial t}=F(u,D_xu,D^2_xu,\cdots),
$$
a stationary point is a solution of the equation
$F(u(x),D_xu(x),D^2_xu(x),\cdots)=0$, at least formally.
The stationary point is a graph on the configuration space.

To extend the concept of the stationary point (or stationary
solution) and establish its existence and the decomposition of
stable and unstable manifolds on the tangent space of the stationary
point to random dynamical systems (RDS) is a basic problem for RDS
(\cite{[Ar]}). In recent years, the stable and unstable manifolds theorem
has been established for finite-dimensional stochastic differential
equations (\cite{[MS3]}) and stochastic systems with memory (see \cite{[MS1]},
\cite{[MS2]}); and for the infinite-dimensional case, the stable manifold
theorem for semilinear stochastic evolution equations and stochastic
partial differential equations (SPDEs) was proved recently by
Mohammed, Zhang and Zhao (\cite{[MZZ]}); invariant manifolds for SPDE and
smooth stable and unstable manifolds for stochastic evolution
equations with one dimensional linear noise were studied by Duan, Lu
and Schmalfuss (\cite{[DLS1]},\cite{[DLS2]}).

To define the pathwise stationary solution of
RDS, let $(\Omega, {\cal F}, P)$ be a complete probability space,
$\theta:\mathbf{\mathbf{R}}\times\Omega\rightarrow\Omega$ be a group
of $P$-preserving ergodic transformations on $(\Omega, {\cal F},
P)$, $H$ be a separable Hilbert space with norm $|\cdot|$ and Borel
$\sigma-$algebra ${\cal B}(H)$. \vskip5pt

 {\bf Definition 1.1} (c.f. \cite{[MZZ]}) A $C^k$ perfect cocyle $(U,\theta)$
on $H$ is a
$(\mathcal{B}(\mathbf{R}^+)\otimes\mathcal{B}(H)\otimes\mathcal{F},\mathcal{B}(H))$-measurable
random field $U:{\mathbf{R}}^+\times H\times\Omega\rightarrow H$
with the following properties:

(i) For each $\omega\in \Omega$, the map ${\mathbf{R}}^+\times H\ni
(t,\xi)\mapsto U(t,\xi,\omega)\in H$ is continuous; for fixed $(t,\omega)\in
{\mathbf{R}}^+\times\Omega$, the map $H\ni \xi\mapsto U(t,\xi,\omega)\in H$ is
$C^k$.

(ii) $U(t_1+t_2,\cdot,\omega)=U(t_2,\cdot,\theta(t_1,\omega))\circ
U(t_1,\cdot,\omega)$ for all $t_1,t_2\in \mathbf{R}^+$, all $\omega\in \Omega$.

(iii) $U(0,\xi,\omega)=\xi$ for all $\xi\in H, \omega\in \Omega$.

\vskip5pt

 {\bf Definition 1.2} (c.f. \cite{[MZZ]}) An $\cal F$ measurable random variable $Y:
\Omega\rightarrow H$ is said to be a stationary random point for the
cocycle $(U,\theta)$ if it satisfies the following identity:
$$
 U(t,Y(\omega),\omega)=Y(\theta(t,\omega))
$$
for all $(t,\omega)\in {\mathbf{R}}^+\times\Omega$. \vskip5pt

 As
pointed out by Mohammed, Zhang and Zhao in \cite{[MZZ]}, this
concept essentially gives a useful realization of the idea of an
invariant measure for the stochastic dynamical system, and allows
us to analyze the local almost sure or generic properties of the
stochastic semiflow in a neighborhood of the stationary point. It
was pointed out by Q. Zhang and Zhao in \cite{[ZZ]}, this
``one-force, one-solution" setting is a natural extension of
equilibrium in deterministic systems to stochastic counterparts.
Unlike the usual serch for invariant measures, it describes the pathwise
invariance of the stationary solution over time along the
measurable and $P$-preserving transformation $\theta_t$:
$\Omega\longrightarrow\Omega$. Only in some special cases, an
invariant measure can be realized as a stationary solution unless one considers
an extended probability space (\cite{[Ar]}). On the other hand, a stationary solution always
generates  an invariant measure. Therefore pathwise stationary solutions
reveal more accurate information about random dynamical systems.    We would like
to point out that unlike deterministic cases, in
general, the stationary solutions of the stochastic systems can
not be given explicitly. Up to now, there has  not been a method
which can be applicable to SPDEs with great generalities. It was
even unthinkable to represent them as solutions of certain
differential or functional equations in general. The existence and
construction of the stationary solutions for SPDEs is a subtle
problem of great importance.

In this article, we will focus on the existence, uniqueness and representation
of the stationary solution of the stochastic Burgers equation (SBE)
with $L^2[0,1]$ valued white noise as follows:
\begin{equation}\label{zhaoa11}
  \left\{
     \begin{array}{ll}
     du(t,x)=(\nu\Delta u(t,x)+\frac{1}{2}\frac{\partial}{\partial
        x}[u(t,x)^2])dt+dW(t), t\geq 0; x\in (0,1)\\
      u(t,1)=u(t,0)=0\\
      u(0,x)=u_0(x).
     \end{array}            \right.
  \label{burgers u}
\end{equation}
This equation has been studied intensively in the literature in the
last ten years (see \cite{[DPZ2]}, \cite{[DP]}, \cite{[EKMS]}, \cite{[MZZ]}, \cite{[S1]}, \cite{[S2]}
and the references therein), because of the interests stemmed from
physics and mathematics. It is a more realistic simple model for turbulence
comparing to deterministic Burgers equations. The latter do not display any chaotic
phenomena as all solutions converge to a unique stationary solution (one graph
on the configuration space). But for the stochastic Burgers equation, we know from
definition that actually the stationary solution is a random moving graphs on the configuration
space (so infinitely many graphs).
In \cite{[S1]}, \cite{[S2]}, Sinai
established the existence and uniqueness of stationary strong
solution of Burgers equations perturbed by periodic forcing or
random forcing. His main tools are Hopf-Cole transformation and
Ito's lemma, hence he required that the noise term which is a
Brownian white noise in time has continuous 3rd-order derivative in spatial
variable. Moreover, he discussed the stationary solution in the
views of the statistical physics. Da Prato and Zabczyk studied the
ergodicity of the SBE in \cite{[DPZ2]} and \cite{[DP]}. In \cite{[MZZ]}, Mohammed, Zhang
and Zhao proved the $C^1$ cocycle property of (\ref{burgers u}). In
 \cite{[EKMS]}, E, Khanin, Mazel and Sinai studied the pathwise stationary solution
 of the stochastic
inviscid Burgers equation with periodic boundary condition.

The main results in this article is to prove under the condition A (large
viscosity condition) which will be made precise in the beginning of
Section 3, there is a unique stationary solution $u$ or stationary
random point $Y(\omega)$ satisfying, for any $t\geq 0$ and a.e. $\omega\in \Omega$,
\begin{equation}
 u(t,\omega)=u(t,Y(\omega), \omega)=Y(\theta(t,\omega)),\label{1.0}
\end{equation}
and
\begin{equation}
  u(t)=\frac{1}{2} \int_{-\infty}^t T_\nu(t-s)\frac{\partial
  u^2(s)}{\partial x}ds+\int_{-\infty}^t T_\nu(t-s)dW_s, \label{1.1}
\end{equation}
or
\begin{equation}
  Y(\omega)=\frac{1}{2}\int_{-\infty}^0T_\nu(-s)\frac{\partial
  (Y(\theta(s,\omega))^2}{\partial
  x}ds+\int_{-\infty}^0T_\nu(-s)dW_s(\omega).
  \label{1.2}
\end{equation}
Moreover, $Y(\omega)$ is hyperbolic and whole $L^2[0,1]$ is its
stable manifold.

Here, we would like to stress the crucial role of the equation
(\ref{1.1}) or (\ref{1.2}) and any stationary solution of the SBE is
given by equation (\ref{1.1}) or (\ref{1.2}). We can only construct
the stationary point under the large viscosity condition. The problem 
remains open without this condition. 
In fact, for SPDE, the essential difficulty of constructing the stationary
point is due to the fact that the equation is
non-autonomous for a.e. $\omega$. But the cocycle
property makes it possible to construct the stationary solution although it is
difficult in general.

For Burgers equation (\ref{zhaoa11}), the
stationary solution should satisfy equation (\ref{1.1}), since it is
easy to see from (\ref{1.1}) that
\begin{equation}
 u(t)=T_\nu(t-t_1)u(t_1)+\frac{1}{2} \int_{t_1}^t T_\nu(t-s)\frac{\partial
  u^2(s)}{\partial x}ds+\int_{t_1}^t T_\nu(t-s)dW_s,
\end{equation}
holds for any $t_1<t$. For other integral equations and backward
doubly stochastic differential equations on infinite horizon, see
 \cite{[MZZ]} and \cite{[ZZ]} respectively.

To construct a solution of (\ref{1.1}), we adopt the pull-back
procedure used by Flandoli and Schmalfuss in \cite{[FS]} in the
construction of the random attractors for 3D Navier-Stokes equations
with non-regular force, and used by Mattingly in \cite{[Ma]} in the
construction of the unique invariant measure of 2D stochastic
Navier-Stokes equation with large viscosity. The main and basic
tools in \cite{[FS]} and \cite{[Ma]} are the so-called Ladyzhenskaya's
inequality (see \cite{[R]} p244) and the weak (or weak$^*$) compactness
method. However, compactness method can not ensure the uniqueness of
the solution of (\ref{1.1}). If uniqueness of the solution of
(\ref{1.1}) does not hold, we cannot prove the stationary point
satisfying (\ref{1.0}). It seems that the uniqueness of the solution
(\ref{1.1}) is a crucial technical condition to obtain the
stationarity of $Y(\omega)$. Therefore, as
in \cite{[Ma]}, we have to require the large
viscosity to ensure the uniqueness of the solution of (\ref{1.1})  and its large time stability.

In fact, even in the case of the deterministic 2D Navier-Stokes
equation with the given exterior force independent of time $t$,
the uniqueness of stationary solution (also called steady-state
solution) is only obtained under the condition of large viscosity (see
Section III.3 in \cite{[So]} or theorem 1.3 in Chapter 2 in \cite{[T2]} for
details, and our condition A is similar to those given in
 \cite{[So]}, \cite{[T2]}).


Let us define some notations used in the rest parts of this article:-
\begin{eqnarray*} L^2[0,1]&\equiv& \{f: [0,1]\to R^1|\
f(0)=f(1)=0,\
\int_0^1f^2(x)dx<\infty\},\\
 H^1_0&\equiv&\{f\in L^2[0,1]|\int_0^1 (\p
f(x))^2dx<\infty\},\\
|f|_p&\equiv &\Big(\int_0^1|f(x)|^pdx\Big)^{\frac{1}{p}}.
\end{eqnarray*}
 In Section 2, we will prove the existence and
uniqueness of the mild solution of (\ref{burgers u}). Here, we
emphasize that the solution of (\ref{burgers u}) under the
$L^2[0,1]$ valued noise is in $L^2([0,T],H^1_0)$. So, combining this property
and the proof of the $C^1$ perfect cocycle in \cite{[MZZ]}, we
have that the mild solution of (\ref{burgers u}) generates a $C^1$
perfect cocycle in $H^1_0[0,1]$. As far as we know, it is not known
that the solution of the 2D stochastic Navier-Stokes equation
generates a $C^1$ perfect cocycle. In this paper we only consider
stochastic Burgers equations. We hope to consider stochastic
Navier-Stokes equations in future publications.

In the section 3, we will construct the stationary solution of SBE.
Some of estimates similar to those of Mattingly (\cite{[Ma]}) are needed.
Although mattingly's estimate is a key
tool for our proof, our idea is different from Mattingly's and equation
(\ref{1.1}) or (\ref{1.2}) is new and plays important roles in our
construction. We also proved the perfection version of (\ref{1.0}) i.e. (\ref{1.0})
holds for all $t$ and a.e. $\omega\in \Omega$. This plays an 
important role in our construction and does
not seem obvious from Mattingly's approach.
In fact, we prove the identity (\ref{1.0}) by using
(\ref{1.1}) and uniqueness of the mild solution of (\ref{burgers
u}). Section 4 is an appendix, in which we will present some
estimates needed in the previous sections.

\section{The weak solution, mild solution and perfect cocycle}

It is well known that $\{e_k(x)=\sqrt{\frac{2}{\pi}}\sin k\pi x,\
x\in[0,1],\ \ k=1,2,\cdots\}$ is a normal orthogonal basis of
$L^2[0,1]$. Let $W(t)$ be a $L^2[0,1]$-valued Brownian motion
defined on the canonical filtered Wiener space $(\Omega,{\cal
F},({\cal F}_t)_{t\in {\bf R}}, P)$. That is, $\Omega$ is the
space of all continuous path $\omega$: ${\bf R}\rightarrow
L^2[0,1]$ such that $\omega(0)=0$ with the compact open topology,
$\cal F$ is its Borel $\sigma$-field,   $P$ is the Wiener measure
on $\Omega$. The Brownian motion is given by:
$$
  W(t,\omega)=\omega(t), \ \ \omega\in \Omega,  \ t\in {\bf R},
$$
and may be represented by $W(t)=\sum_{k=1}^\infty\sigma_k e_k
B_k(t)$, where $\sum\limits_{k=1}^{\infty}\sigma _k^2<\infty$, $\{B_k(t),\ t\in(-\infty,\infty) \ \
k=1,2,\cdots\}$ are mutually independent real valued Brownian
motions on $(\Omega,{\cal F}, P)$,  and
$$
  {\cal F}_t=\sigma\{\omega(s),s\leq t\}\vee{\cal N}, \ \  \
  {\rm for \  any} \ \ t\in {\bf R}
$$
where $\cal N$ are the null sets of $\cal F$ (see \cite{[DPZ1]} p86-87
and \cite{[Ar]} p91).

Throughout this article, we denote by $\theta$: ${\bf
R}\times\Omega\rightarrow\Omega$ the standard $P$-preserving
ergodic Wiener shift on $\Omega$
$$
\theta(t,\omega)(s)\equiv \omega(t+s)-\omega(t), \ \ t,s \in {\bf R
}.
$$
Hence $(W, \theta)$ is an helix:
$$
W(t_1+t_2,\omega)-W(t_1,\omega)=W(t_2,\theta(t_1,\omega)),\ \
t_1,t_2\in {\bf R }, \ \omega\in \Omega,
$$
and
\begin{equation}
\theta^{-1}(u){\cal F}_t={\cal F}_{t+u}, \ \ \  \forall\ t,u\in
{\bf R}. \label{theta}
\end{equation}


Now, we denote $\displaystyle \w(t)=\int_0^t T_\nu(t-s)dW(s)$, and
$\displaystyle \w(t_0,t)=\int_{t_0}^t T_\nu(t-s)dW(s)$, the
so-called stochastic convolution with respect to the semigroup
$T_\nu(t)=e^{-t\nu\Delta}$. Some properties of the stochastic
convolution $\w$ are presented in Appendix A. In particular, theorem
A.2 shows that $\w\in L^2([0,t], H^1_0)\ \ a.s.$.

Now, consider the following PDE with random coefficients
$$
v(t)=T_\nu(t)v_0+\frac{1}{2}\int_0^t T_\nu(t-s)[\frac{\partial}{\partial
x}(v(s)+\w (s))^2]ds,
$$
which can be considered as the mild form of the random Burgers
equation
\begin{equation}
  \left\{
     \begin{array}{ll}
      \frac{dv(t)}{dt}=\nu\Delta v(t)+\frac{1}{2}\frac{\partial}{\partial
        x}[v(t)+\w(t)]^2,\\
      v(0)=v_0(x).
     \end{array}            \right.
  \label{burgers v}
\end{equation}

Let's consider the ODE used in Galerkin approximation
\begin{equation}
  \left\{
    \begin{array}{ll}
       \frac{dv_n(t)}{dt}=\nu\Delta v_n(t)+\frac{1}{2}P_n\frac{\partial}{\partial
        x}[v_n(t)+\w(t)]^2,\\
        v_n(t_0)=P_n(u_{t_0}),
    \end{array}          \right.
  \label{A1}
\end{equation}
where $P_n$ is the projection operator in $L^2[0,1]$ onto the space
spanned by $\{e_1,e_2,\cdots, e_n\}$, and $\w(t_0)=0$.

\medskip

{\bf Lemma 2.1}: {\it The global solution to (\ref{A1}) exists and is unique.}
\medskip

{\bf Proof:} It's easy to show the local existence and uniqueness by the
standard fixed point theorem, and then
$$
\langle v_n, \frac{dv_n(t)}{dt}\rangle=\langle v_n(t),\Delta
v_n(t)\rangle+\langle
v_n(t),\frac{1}{2}P_n[\frac{\partial}{\partial
x}[v_n(t)+\w(t)]^2]\rangle.
$$
By the Dirichlet boundary conditions on $[0,1]$ and integration by
parts
\begin{eqnarray}
   \frac{1}{2}\frac{d}{dt}|v_n(t)|^2_2+|\nabla v_n(t)|^2_2
   & = &\frac{1}{2}\langle v_n(t),\frac{\partial}{\partial
   x}[v_n(t)+\w(t)]^2\rangle \label{A2}\nonumber \\
   & = &-\int_0^1v_n(t,x)\w(t,x)\frac{\partial}{\partial
   x}v_n(t,x)dx\nonumber\\
   &&
   -\frac{1}{2}\int_0^1 \w^2(t,x)\frac{\partial}{\partial
   x} v_n(t,x)dx \nonumber\\
   & = &{\rm I}+{\rm II}.
\end{eqnarray}
By the H\"{o}lder inequality and Sobolev embedding theoem, we have
$$
|{\rm II}|=\Big|\int_0^1 \w^2(t,x)\frac{\partial}{\partial
   x}v_n(t,x)dx\Big|\leq |\w(t)|^2_4\Big|\frac{\partial}{\partial
   x}v(t)\Big|_2;
$$
and
\begin{eqnarray*}
 |\rm {I}|^2& = &\Big|\int_0^1v_n(t,x)\w(t,x)\p v_n(t,x)dx\Big|^2\\
      &\leq&|v_n(t)|_4|\w(t)|_4\Big|\p v_n(t,\cdot)\Big|_{L^2[0,1]}\\
      &\leq&\gamma^2|\nabla v_n(t)|^{1\over 2}_2|v_n(t)|^{1\over
      2}_2|\w(t)|_4|\nabla v_n(t)|_2\\
      & = &\gamma^2|\nabla v_n(t)|^{3\over 2}_2|v_n(t)|^{1\over
      2}_2|\w(t)|_4.
\end{eqnarray*}
Hence, due to Young's inequality,
\begin{eqnarray*}
 &    &   {1\over 2}{d\over dt}|v_n(t)|^2_2+|\nabla
 v_n(t)|^2_2\\
 &\leq& \varepsilon|\nabla
 v_n(t)|^2_2+C_1(\varepsilon)|v_n(t)|^2_2|\w(t)|^4_4
 +C_2(\varepsilon)|\w(t)|^4_4.
\end{eqnarray*}
Taking $\varepsilon={1\over 2}$, we obtain
$$
{d\over dt}|v_n(t)|^2_2+|\nabla v_n(t)|^2_2\leq
C|v_n(t)|^2_2|\w(t)|^4_4+|\w(t)|^4_4
$$
for a constant $C>0$. By the Gronwall inequality,
\begin{eqnarray}
 |v_n(t)|^2_2 &\leq & |v_n(t_0)|^2_2\exp\Big\{\int_{t_0}^t
 C|\w(s)|^4_4ds\Big\}\nonumber
 \\
 &&
 +\int_{t_0}^t\exp\big\{\int_r^t
 C|\w(s)|^4_4ds\big\}\,|\w(r)|^4_4dr\label{A2'},
\end{eqnarray}
and
\begin{equation}
\int_r^t|\nabla
v_n(s)|^2_2ds\leq|v_n(r)|^2_2+C\int_r^t(|v_n(s)|^2_2|\w(s)|^4_4
+|\w(s)|^4_4)ds\label{A3}.
\end{equation}
Using theorem 5.20 in \cite{[DPZ1]} (p141) that says that $\w(\cdot)\in
C([0,T];C[0,1])$ and from (\ref{A2'}),(\ref{A3}),we have
\begin{equation}
  \sup_{t\in[0,T]}|v_n(t)|_2\leq K_1(\omega) \ \ \mbox{\rm and}\ \  \int_0^T|\nabla
  v_n(s)|^2_2ds\leq K_2(\omega)\ \ \mbox{\rm uniformly in}\ \ n.\label{A4}
\end{equation}
The  first inequality implies that $|v_n(t)|_2$ is finite for all
$t>0$. Hence, by Lemma 2.4 in \cite{[R]} (p48), we have the global solution
to (\ref{A1}).\hfill $\Box$
\medskip

Inequalities in (\ref{A4}) imply that $v_n$ is bounded in $L^\infty(0,T; L^2[0,1])$ and $L^2(0,T; H^1_0)$
uniformly in
$n$. These
uniform bounds allow us to use the Alaoglu compactness theorem to
find a subsequence which we shall denote it by $\{v_n\}$ such that
\begin{equation}
 v_n\stackrel{*}\rightarrow v\ \ {\rm in}\ \ L^\infty(0,T; L^2[0,1])
\end{equation}
{\it i.e.} $v_n$ weak$^*$ converge to $v$. We can extract a further
subsequence, still denote by  $\{v_n\}$, such that
\begin{equation}
 v_n\rightarrow v\ \ {\rm in}\ \ L^2(0,T; H^1_0)
\end{equation}
with
\begin{equation}
 v\in  L^\infty(0,T; L^2[0,1])\cap L^2(0,T; H^1_0).\label{2.6}
\end{equation}
From the discussion of Theorem 6.1 in \cite{[R]}, we know that the
Laplacian operator $\Delta$ is bounded from $H^1_0$ to $H^{1*}_0$
in the sense defined in Theorem 6.1 in \cite{[R]}. Since $v_n$ is bounded
uniformly in $L^2(0,T; H^1_0)$, $\nu \Delta v_n$ is bounded
uniformly in $L^2(0,T; H^{1*}_0)$.

Now, let us prove that $P_n[\p v_n^2]$ are uniformly bounded in
$L^2(0,T; H^{1*}_0)$.  First it is easy to show $|\p
v_n^2|_{H^{1*}_0}\leq C|v_n|_2|\nabla v_n|_2$. So Lemma 7.5
in \cite{[R]} implies
$\|P_nB(v,v)\|_{L^2(0,T;H^{1*}_0)}\leq\|B(v,v)\|_{L^2(0,T;H^{1*}_0)}$.
Therefore,
\begin{eqnarray}
  \|P_nB(v_n,v_n)\|_{L^2(0,T;H^{1*}_0)}&\leq&\int_0^T|B(v_n,v_n)|_{H^{1*}_0}^2ds\nonumber\\
  &\leq&C\int_0^T|v_n(s)|^2|\nabla v_n(s)|^2ds \label{A3'}\\
  &\leq&C\|v_n(s)\|^2_{L^\infty(0,T;L^2)}\|v_n\|^2_{L^2(0,T;H^{1*}_0)}.
  \nonumber
\end{eqnarray}
\vskip5pt

{\bf Lemma 2.2}. {\it  The random perturbed
$P_nB(v_n+W_\Delta,v_n+W_\Delta)$ is uniformly bounded in
$L^2(0,T;H^{1*}_0)$ a.s..}
\medskip

{\bf Proof:} As $v_n$ is uniformly bounded in $L^\infty(0,T;
L^2[0,1])$ and $L^2(0,T;H^{1*}_0)$, so does $v_n+W_\Delta$. By
(\ref{A3'}), we finish the lemma. \hfill $\Box$

\vskip5pt


{\bf Lemma 2.3:} {\it There exists a subsequence $v_n$ such that ${{dv_n}\over
{dt}}\stackrel{*}\rightarrow {{dv}\over {dt}}$ in $L^2(0,T;H^{1*}_0)$}.
\medskip

{\bf Proof:} Lemma 2.1 implies ${dv_n}\over {dt}$ is uniformly
bounded in $L^2(0,T;H^{1*})$. We can extract a further subsequence
(relabelling again). Using the same argument as in \cite{[R]} (p203-p204),
we have
$$
{{dv_n}\over {dt}}\stackrel{*}\rightarrow {{dv}\over {dt}}\ \ {\rm in}\ \
L^2(0,T;H^{1*}_0).
$$
\hfill$\Box$


{\bf Lemma 2.4:} {\it There exists a subsequence $\{v_n\}$ such that
}
$$
P_n\Big[{{dv_n}\over {dt}}\Big]\stackrel{*}\rightarrow \Big[{{dv}\over
{dt}}\Big]\ \ {\rm in}\ \ L^2(0,T;H^{1*}_0).
$$
\medskip

{\bf Proof:} {Step 1}: Since $H^1_0\subset\subset
L^2[0,1]$, by Theorem 8.1 in \cite{[R]}, there is a subsequence $\{v_n\}$
(after relabelling) that converges to $v$ strongly in $L^2(0,T;
L^2[0,1])$.

{Step 2}: Let $K=\{\sum_{i=1}^k\alpha_i(t)\phi_i;\
\phi_i\in H^1_0, \alpha_i\in C[0,T], k\in N \}$. We will show
for any $\phi \in K$,
\begin{equation}
\int_0^T\int_0^1 v_n(\p v_n)\phi \,dxdt\rightarrow\int_0^T\int_0^1 v(\p
v)\phi\, dxdt. \label{A9}
\end{equation}
In fact,
\begin{eqnarray}\label{zhao1}
   &&\Big|\int_0^T\int_0^1 v_n(\p v_n)\phi \,dxdt-\int_0^T\int_0^1 v(\p
v)\phi\, dxdt\Big|\nonumber\\
   &\leq &\sum_{i=1}^k\Big|\int_0^T\alpha_i\left [\int_0^1v_n(\p v_n)\phi_i \,dx-\int_0^1 v(\p
v)\phi_i\, dx\right ]dt\Big|\nonumber\\
   &\leq
   &{{1}\over{2}}\sum_{i=1}^k\Big|\int_0^T\alpha_i\Big[\int _0^1(\p\phi_i)v_n^2dx
   -\int_0^1(\p\phi_i)v^2dx\Big]dt\Big|\nonumber\\
   &\leq&{{1}\over{2}}\sup_{i;t\in[0,T]}|\alpha_i(t)|\sum_{i=1}^k\int_0^T\int_0^1\Big|\p
   \phi_i\Big||v_n^2-v^2|dxdt.
   \end{eqnarray}
  Then by the Sobolev embedding theorem and Cauchy-Schwartz inequality,
  \begin{eqnarray}\label{zhao2}
   &&\int_0^T\int_0^1\Big|\p
   \phi_i\Big||v_n^2-v^2|dxdt\nonumber\\
     &\leq&C\int_0^T|v_n+v|_{H^1_0}\int_0^1\Big|\p
   \phi_i\Big||v_n-v|dxdt\nonumber\\
   &\leq&C\int_0^T|v_n+v|_{H^1_0}[\int_0^1\Big|\p
   \phi_i\Big|^2dx]^{1\over 2}|v_n-v|_2dt\nonumber\\
   &\leq&C[\int_0^1\Big|\p
   \phi_i\Big|^2dx]^{1\over
   2}\left [\int_0^T|v_n+v|_{H^1_0}^2dt\right]^{1\over 2}
   \left [\int^T_0 |v_n-v|^2_2dt\right ]^{1\over 2}.
\end{eqnarray}
It is easy to see from (\ref{zhao1}) and (\ref{zhao2}) that
\begin{eqnarray*}
\Big|\int_0^T\int_0^1 v_n(\p v_n)\phi \,dxdt-\int_0^T\int_0^1 v(\p
v)\phi\, dxdt\Big|\leq C \left [\int^T_0 |v_n-v|^2_2dt\right
]^{1\over 2}.
\end{eqnarray*}
Here $C$ is a generic positive constant.  Since $v_n$ is uniformly bounded in $L^2(0,T;H^1_0)$,
the convergence in (\ref{A9}) follows.

{Step 3}: It is obvious that
\begin{eqnarray*}
  &&\Big|\int_0^T\int_0^1 P_n v_n(\p v_n)\phi \,dxdt-\int_0^T\int_0^1 v(\p
v)\phi\, dxdt\Big|\\
  &\leq&\Big|\int_0^T\sum_{i=1}^k \Big(\int_0^1(\p v_n^2)P_n\phi_idx\alpha_i(t)dt-\int_0^T\int_0^1(\p
  v_n^2)\phi_idx\alpha_i(t)dt\Big)\Big|\\
  &  &+\Big|\int_0^T\sum_{i=1}^k\Big(\int_0^1(\p v_n^2)\phi_idx\alpha_i(t)dt-\int_0^1(\p v^2)\phi_idx\alpha_i(t)dt\Big)\Big|\\
  &= &{\rm I}+{\rm II}.
\end{eqnarray*}
By Step 2, we know II$\rightarrow0$ as $n\to \infty$. To see
I$\rightarrow0$ as $n\to \infty$, note
\begin{eqnarray*}
&&\Big|\int _0^T\int_0^1\Big[\p\Big(
P_n\phi_i\alpha_i(t)-\phi_i\alpha_i(t)\Big)\Big]^2dxdt\Big|\\
&\leq&
\int_0^T\alpha_i^2(t)dt\int_0^1\Big|\p
(P_n\phi_i-\phi_i)\Big|^2dx\to 0,
\end{eqnarray*}
noticing that Lemma 7.5 in \cite{[R]} implies $P_n\phi\rightarrow\phi$ in
$L^2(0,T;H^0_1)$. Moreover
 $${\rm I}\leq\|v_n(\p
v_n)\|_{L^2(0,T;H^{1*}_0)}\|P_n\phi-\phi\|_{L^2(0,T;H^0_1)}.
$$
 Since
(\ref{A3'}) means that $v_n(\p v_n)$ is uniformly bounded in
$L^2(0,T;H^{1*}_0)$, we have ${\rm I}\rightarrow0$ as $n\to
\infty$.

{Step 4}: It is well known that $K$ is dense in
$L^2(0,T;H^1_0)$. Therefore, Step 2 and Step 3 show
that for any for any $\phi\in L^2(0,T;H^1_0)$,
\begin{equation}
  \int_0^T\int_0^1 P_n v_n(\p v_n)\phi
  \,dxdt\rightarrow\int_0^T\int_0^1 v(\p v)\phi\, dxdt.
\end{equation}
This proves the lemma. \hfill $\Box$
\medskip

Using a similar argument as in \cite{[R]} (p249), we have $v_n(0)=P_n
u_0\rightarrow v_0=u(0)$. Therefore, using the above lemmas, it is
easy to deduce the following theorem.


\vskip5pt

{\bf Theorem 2.5} {\it There exists a solution $v\in
L^2(0,T;H^1_0)\cap L^\infty(0,T; L^2[0,1])$
 that satisfies
\begin{equation}\label{zhao3}
  {dv\over dt}=\nu\Delta v+ \frac{1}{2}{\partial(v+\w)^2\over \partial x}
\end{equation}
as an equation in $L^2(0,T;H^{1*}_0)$. }

\vskip5pt



{\bf Theorem 2.6} {\it There exists a unique mild solution $u\in
L^2(0,T;H^1_0)\cap C(0,T; L^2[0,1]) $ that satisfies
\begin{equation}
u(t)=T_\nu(t)u_0+{1\over 2}\int_0^tT_\nu(t-s)[\p
u^2(s)]ds+\int_0^tT_\nu(t-s)dW(s), \label{solution}
\end{equation}
where $u_0\in L^2[0,1]$.}

\vskip5pt

Before we prove theorem 2.6, we need some lemmas. Denote $p(t,x,y)$
the heat kernel of $\nu\Delta$ with the Dirichlet boundary
conditions on $[0,1]$. Note here $p(\cdot,\cdot,y)\not\in
L^2(0,T;H^1_0)$.

Let $f\in C_0^\infty([0,1])$ and $K(t,y)=\int_0^1 p(t,x,y)f(x)dx$,
then $K(t,y)\in C^\infty((0,T]\times[0,1])$ and
$K(t,y)=\int_0^1p(t,y,x)f(x)dx$. Obviously, the energy inequality
implies
\begin{equation}
  K(t,y) \in L^2(0,T;H^1_0).\label{A10}
\end{equation}


\medskip

{\bf Lemma 2.7} {\it Let $v\in L^2(0,T;H^1_0)\cap L^\infty(0,T;
L^2[0,1])$ be the solution of equation (\ref{zhao3}) in
$L^2(0,T;H^{1*}_0)$, then $\int_0^t\int_0^1 p(t-s,x,y)
\frac{\partial}{\partial y}(v+\w)^2(s,y)dyds\in L^2[0,1]$. }

\medskip

{\bf Proof:} For every $f\in L^2[0,1]$, by Cauchy-Schwartz
inequality
\begin{eqnarray*}
 & &\int_{[0,T]\times[0,1]\times[0,1]}|p(t-s,x,y)f(x)\frac{\partial}{\partial
 y}(v+\w)^2(s,y)|dxdyds\\
 &\leq&2\Big[\int_{[0,T]\times[0,1]\times[0,1]}p(t-s,x,y)f^2(x)(v+\w)^2(s,y)dxdyds\Big]^{1\over
 2}\\
 & &\times \Big[\int_{[0,T]\times[0,1]\times[0,1]}p(t-s,x,y)\Big(\frac{\partial}{\partial
 y}(v+\w)\Big)^2(s,y)dxdyds\Big]^{1\over 2}\\
 &\leq&2C\Big(\int_0^t{1\over {\sqrt{t-s}}}ds\Big)^{1\over 2}
 \Big(\int_0^1 f^2(x)dx\Big)^{1\over
 2}\Big[\|v+\w\|_{L^\infty(0,T;L^2[0,1])}\Big]^{1\over
 2}\\
 & &\times \Big[\int^t_0\int^1_0 \left(\int_0^1p(t-s,x,y)dx\right)
 \Big(\frac{\partial}{\partial
 y}(v+\w)\Big)^2dyds\Big]^{1\over 2}\\
 & \leq&C(t)\|v+\w\|^{1\over
 2}_{L^\infty(0,T;L^2[0,1])}\cdot \|v+\w\|_{L^2(0,T;H^1_0)}\cdot\|f\|^2_{L^2([0,1])}.
\end{eqnarray*}
This means that the linear functional ${\cal L}: f\rightarrow
\int_0^1 f(x)\int_0^t\int _0^1
p(t-s,x,y)\Big(\frac{\partial}{\partial
 y}(v+\w)(s,y)\Big)^2dydsdx$ is bounded in $L^2[0,1]$. Then by
the Rieze representation theorem, we have $$\int_0^t\int_0^1
p(t-s,x,y)\Big(\frac{\partial}{\partial
 y}(v+\w)(s,y)\Big)^2dyds\in L^2[0,1].$$
\hfill $\Box$
\medskip



{\bf Lemma 2.8} {\it The following equality holds in $L^2[0,1]$}
\begin{equation}
v(t)=T_\nu (t)v(0)+{1\over 2}\int_0^t T_\nu(t-s)\p (v+\w)^2ds
\end{equation}

{\bf Proof:} From (\ref{A10}) and Theorem 2.5, we know
\begin{eqnarray*}
 \int^t_0\int_0^1 K(t-s,x){\partial v(s,x)\over \partial s}dxds&=&-\int_0^t\int_0^1\p
  K(t-s,x)\p v(s,x)dxds\\
  &&+{1\over 2}\int_0^t\int_0^1
  K(t-s,x)\p(v+\w)^2dxds.
\end{eqnarray*}
It follows from the integration by parts formula and the fact that
\begin{equation}
  \left\{
     \begin{array}{ll}
      \frac{\partial}{\partial t} K(t,y)={1\over 2}\Delta K(t,y),\\
      K(t,0)=K(t,1)=0,\\
      K(0,y)=f(y),
     \end{array}            \right.
\end{equation}
then
\begin{eqnarray*}
  & &\int_0^1 f(y)v(t,y)dy\\
  &=&\int_0^1K(t,y)v(0,y)dy+{1\over 2}\int_0^t\int_0^1\int_0^1p(t-s,y,x)f(y)dy\p
    (v+\w)^2dxds\\
    &=&\int_0^1f(y)\int_0^1p(t,y,x)v(0,x)dxdy\\
    &&
    +{1\over
    2}\int_0^1f(y)\int_0^t\int_0^1 p(t-s,y,x)\p
    (v+\w)^2dxdsdy.
\end{eqnarray*}
Due to Lemma 2.7 and $C_0^\infty([0,1])$ is dense in $L^2[0,1]$, we
get the following equality in $L^2[0,1]$,
$$
v(t)=T_\nu(t)v(0)+{1\over 2}\int_0^t T_\nu(t-s)\p (v+\w)^2ds.
$$
\hfill $\Box$


{\bf The proof of theorem 2.6:}  Lemma 2.8 and the definition of
$\w$ imply the existence of a mild solution. Now, we only need to
show that $u\in C(0,T; L^2[0,1])$ and the uniqueness.

By (\ref{2.6}), we know $v\in L^2(0,T; H^1_0)$ and by Lemma 2.3 we
know ${{dv}\over {dt}}$ in $L^2(0,T;H^{1*}_0)$. Thus Theorem 7.2
in \cite{[R]} implies that $v\in C(0,T; L^2[0,1])$. At the same
time, by Theorem 5.20 in \cite{[DPZ1]}, $\w \in C([0,T];C[0,1])$.
It then follows that $u\in C(0,T; L^2[0,1])$.

Let $u_1$ and $u_2$ be two solutions of (\ref{solution}) with the same initial
data, then using Proposition A.3 in the Appendix, we have
\begin{eqnarray*}
 & &|u_1(t)-u_2(t)|_2^2
 =\int_0^1\Big(u_1(t,x)-u_2(t,x)\Big)^2dx\\
 &=&\int_0^1[\int_0^t T_\nu(t-s)\Big(\frac{\partial u_1^2}{\partial y}-
  \frac{\partial u_2^2}{\partial y}\Big)ds]^2dx\\
 &\leq& C\int_0^1dx\Big(\frac{1}{2}\int_0^t
   ds\int_0^1\frac{1}{\sqrt{t-s}}\frac{c_1}{\sqrt{t-s}}e^{-\frac{(x-y)^2}{2c_2(t-s)}}
   |u_1^2(s,y)-u_2^2(s,y)|dy\Big)^2\\
 &\leq& C\int_0^1\int_0^t\frac{1}{(t-s)^{\frac{3}{4}}}ds\int_0^t\frac{1}{(t-s)^{\frac{1}{4}}}
 \\
 &&\hskip2cm \Big(\int_0^1\frac{c_1}{\sqrt{t-s}}e^{-\frac{(x-y)^2}{2c_2(t-s)}}|u_1^2(s,y)-u_2^2(s,y)|dy\Big)^2
 dsdx\\
 &\leq& C\int_0^1\int_0^t\frac{1}{(t-s)^{\frac{1}{4}}}\int_0^1\frac{c_1}{\sqrt{t-s}}e^{-\frac{(x-y)^2}{2c_2(t-s)}}
            (u_1(s,y)-u_2(s,y))^2dy\\
 & &
 \times\int_0^1\frac{c_1}{\sqrt{t-s}}e^{-\frac{(x-y)^2}{2c_2(t-s)}}(u_1(s,y)+u_2(s,y))^2dydsdx\\
 &\leq&C\sup_{0\leq s\leq
 t}|u_1(s)+u_2(s)|^2_2\int_0^t\frac{1}{(t-s)^{\frac{3}{4}}}|u_1(s)-u_2(s)|^2_2ds\\
 &\leq&C(\omega)\int_0^t\frac{1}{(t-s)^{\frac{3}{4}}}|u_1(s)-u_2(s)|^2_2ds,
\end{eqnarray*}
where $C$ is a generic constant that may change from one line to another. We iterate
the above computation,
$$
 |u_1(t)-u_2(t)|_2^2
 \leq
 C(\omega)^2\int_0^t\frac{1}{(t-s)^{\frac{3}{4}}}
 \int_0^s\frac{1}{(s-r)^{\frac{3}{4}}}
 |u_1(r)-u_2(r)|_2^2drds.
$$
Consider now the elementary estimate
$$
\int_r^t\frac{s^\alpha}{(t-s)^\beta(s-r)^\gamma}ds=
\int_0^{t-r}\frac{(s+r)^\alpha}{(t-r-s)^\beta s^\gamma}ds\leq
\frac{C}{(t-r)^{\beta+\gamma-1}}, \ t\geq r>0,
$$
we have
$$
 |u_1(t)-u_2(t)|_2^2\leq
 C\int_0^t\frac{1}{(t-r)^{\frac{1}{2}}}|u_1(r)-u_2(r)|_2^2dr.
$$
Iterating it, and using the  elementary estimate again, we get
$$
|u_1(t)-u_2(t)|_2^2\leq C\int_0^t|u_1(s)-u_2(s)|_2^2ds.
$$
So, the Gronwall inequality implies that $u_1=u_2$ in $L^2[0,1]$. We
complete the uniqueness. \hfill $\Box$

\vskip5pt

Moreover, the mild solution of (\ref{burgers u}) generates a perfect $C^1$
cocycle.
\medskip


{\bf Theorem 2.9} (see also Theorem 1.4.3 \cite{[MZZ]}) {\it
Consider stochastic Burgers equation (\ref{burgers u}). Then
 equation (\ref{burgers u}) has a unique mild solution with a
 $\big({\cal B}({\bf R}^+)\otimes{\cal B}(L^2[0,1])\otimes{\cal F}, {\cal
 B}(L^2[0,1])\big)$ measurable version $u$: ${\bf R}^+\times L^2[0,1]\times\Omega\rightarrow L^2[0,1]$
having the following properties:

(i) For each $\psi\in L^2[0,1], u(\cdot,\psi,\cdot): {\bf
R}^+\times\Omega\rightarrow L^2[0,1]$ is $({\cal F}_t)_{t\geq0}$
adapted;

(ii) $(u, \theta)$ is a $C^1$ perfect cocycle on $L^2[0,1]$;

(iii) For each $(t,\omega)\in (0,\infty)\times\Omega$, the map
$L^2(0,1)\ni \psi\mapsto u(t,\psi,\omega)\in L^2[0,1]$ takes
bounded set into relatively compact sets;

(iv) For each $(t,\psi, \omega)\in (0,\infty)\times
L^2[0,1]\times\Omega$, the Fr\'{e}chet derivative $D
u(t,\psi,\omega)\in L(L^2[0,1])$ is compact. Furthermore, the map
$$
[0,\infty)\times L^2[0,1]\times\Omega\ni(t,\psi, \omega)\mapsto D
u(t,\psi,\omega)\in L(L^2[0,1])
$$
is strong measurable and for each $(t,\omega)\in
(0,\infty)\times\Omega$, the map $L^2[0,1]\ni\psi\mapsto
Du(t,\psi,\omega)\in L^2[0,1]$ takes bounded sets into relative
compact sets;

(v) For any positive $\alpha,\rho$
$$
E\log^+\sup_{0\leq t\leq \alpha, |\psi|_2\leq\rho}
\Big\{|u(t,\psi,\cdot)|_2+\|Du(t,\psi,\cdot)\|_{L(L^2[0,1])}\Big\}<\infty.
$$
} {\bf Proof}: The uniqueness of the solution implies that the
Galerkin approximation sequences (\ref{A1}) $v_n$ converge to $v$
for a.e. $\omega$ in weak$^*$ topology $L^\infty((0,T),L^2[0,1])$.
This shows that $v$ is $\big({\cal B}({\bf R}^+)\otimes{\cal
B}(L^2[0,1])\otimes{\cal F}, {\cal
 B}_w(L^2[0,1])\big)$ measurable, where ${\cal B}_w(L^2[0,1])\big)$ is the
 $\sigma$-algebra generated by the weak topology in $L^2[0,1]$. Because $L^2[0,1]$
 is a separable Hilbert space, it is well-known that $v$ is $\big({\cal
B}({\bf R}^+)\otimes{\cal B}(L^2[0,1])\otimes{\cal F}, {\cal
 B}(L^2[0,1])\big)$ measurable. Since $u=v-W_{\nu\Delta}$ and $W_{\nu\Delta}$
 is $\big({\cal
B}({\bf R}^+)\otimes{\cal B}(L^2[0,1])\otimes{\cal F}, {\cal
 B}(L^2[0,1])\big)$ measurable, we know that $u$ is $\big({\cal
B}({\bf R}^+)\otimes{\cal B}(L^2[0,1])\otimes{\cal F}, {\cal
 B}(L^2[0,1])\big)$.

The proof of (i)--(v) is given in \cite{[MZZ]} in details. \hfill $\Box$


\section{The construction of the stationary solution}
We define $u(t,\omega; t_0, u_0)$ the value of the mild solution of
(\ref{burgers u}) at time $t$ with an initial data $u_0$ at time
$t_0$. We define $\delta_0=\lambda_1\nu-\frac{\varepsilon_0\gamma}{2\nu^2}$,
where $\lambda_1=\pi^2$, the first non-zero eigenvalue of $\Delta$;
$\varepsilon_0\equiv E(|W_1-W_0|^2_2)=\sum_{k=1}^\infty \sigma_k^2$; $\gamma$
is the minimal constant such that the inequality, $\max_{x\in[0,1]}|u(x)|\leq
\gamma \|u\|_{H^1_([0,1])}$ if $u\in H^1([0,1])$, holds.

In this section, we assume
\medskip

{\bf Condition A}:
$\frac{\nu^3}{\varepsilon_0}>\frac{\gamma}{2\lambda_1}$ i.e. $
\delta_0>0$.
\medskip

\noindent
Then we can prove the following lemma.
\medskip

{\bf Lemma 3.1}: {\it Assume that $u(t)$ satisfies the following
equation in $L^2[0,1]$ for any $t\in {\bf R}$
\begin{equation}\label{zhao11}
  u(t)=\frac{1}{2} \int_{-\infty}^t T_\nu(t-s)\frac{\partial
  u^2(s)}{\partial x}ds+\int_{-\infty}^t T_\nu(t-s)dW_s \label{3.1}
\end{equation}
and
\begin{equation}
 \sup_{t\in {\bf R}}E|u(t)|^{2p}_2<\infty, \ \ p\geq 1, \label{3.1c}
\end{equation} then
$u(t)$ is unique.}
\medskip

{\bf Proof}: For any $t_1<t$, by Fubini Theorem, we have
\begin{eqnarray*}
u(t)&=&\frac{1}{2} \int_{-\infty}^{t_1} T_\nu(t-s)\frac{\partial
  u^2(s)}{\partial x}ds+\int_{-\infty}^{t_1} T_\nu(t-s)dW_s\\
  & &\ \ \ \ +\frac{1}{2} \int_{t_1}^t T_\nu(t-s)\frac{\partial
  u^2(s)}{\partial x}ds+\int_{t_1}^t T_\nu(t-s)dW_s\\
  &=&T_\nu(t-t_1)\Big(\frac{1}{2} \int_{-\infty}^{t_1} T_\nu(t_1-s)\frac{\partial
  u^2(s)}{\partial x}ds+\int_{-\infty}^{t_1} T_\nu(t_1-s)dW_s\Big)\\
  & &\ \ \ \ +\frac{1}{2} \int_{t_1}^t T_\nu(t-s)\frac{\partial
  u^2(s)}{\partial x}ds+\int_{t_1}^t T_\nu(t-s)dW_s\\
  &=&T_\nu(t-t_1)u(t_1)+\frac{1}{2} \int_{t_1}^t T_\nu(t-s)\frac{\partial
  u^2(s)}{\partial x}ds+\int_{t_1}^t T_\nu(t-s)dW_s.
\end{eqnarray*}
Therefore, for any $t_1<t\in {\bf R}$, $u(\cdot)$ is the mild solution of the
equation (\ref{burgers u}) with initial data $u(t_1)$.

Now, let's show the uniqueness of the equation (\ref{3.1}). Assume $u(\cdot)$
and $l(\cdot)$ are two solutions of (\ref{3.1}) under condition (\ref{3.1c}), we
have, for any $n\in {\bf Z}^+$, $-n<t$,
$$
u(t)=T_\nu(t+n)u(-n)+\frac{1}{2} \int_{-n}^t T_\nu(t-s)\frac{\partial
  u^2(s)}{\partial x}ds+\int_{-n}^t T_\nu(t-s)dW_s,
$$
$$
l(t)=T_\nu(t+n)l(-n)+\frac{1}{2} \int_{-n}^t T_\nu(t-s)\frac{\partial
  l^2(s)}{\partial x}ds+\int_{-n}^t T_\nu(t-s)dW_s.
$$Since $u$ and $l$ are the mild solutions of equation (\ref{burgers u}) with
initial data $u(-n)$ and $l(-n)$ at time $-n$, by the theorem A.6
in Appendix, it is easy to know, for any $t\in {\bf Z}$ and any
$\varepsilon>0$,  there exist ${\bf Z}$-valued random time
$\overleftarrow{n}$, such that for any $n>\overleftarrow{n}$ with
probability 1,
\begin{equation}
|u(t_1+\tau, t_1-n, u(t_1-n),\omega)-l(t_1+\tau, t_1-n,
l(t_1-n),\omega)|_2^2\leq \varepsilon\delta^2|n|e^{-\delta
(n+\tau)}, \label{3.1d}
\end{equation}
where $\delta \in (0,\delta _0)$ is a constant. So,
\begin{equation}
 |u(t_1+\tau,\omega)-l(t_1+\tau,\omega)|_2^2<\varepsilon, \ \ \mbox{for any}\ \tau>0.\label{3.1e}
\end{equation}
This shows that $u(\cdot)=l(\cdot)$ for $a.e.\omega$. \hfill $\Box$



\medskip

{\bf Lemma 3.2}: {\it Let $u$ satisfy (\ref{3.1}) and (\ref{3.1c}) in Lemma
3.1, then we have, for any $r\in {\bf R}$
\begin{equation}
 u(\cdot, \theta(r,\omega))=u(\cdot+r,\omega)\ \ {\rm for \ a.e.}\ \omega. \label{3.2}
\end{equation}
}

 {\bf Proof}: Because $\theta$ is a $P-$preserving on the probability space
$(\Omega,{\cal F}, P)$, for any $r$, we have that the equation
$$
 u(t,\theta(r,\omega))=\frac{1}{2} \int_{-\infty}^t T_\nu(t-s)\frac{\partial
  u^2(s,\theta(r,\omega))}{\partial x}ds+\int_{-\infty}^t T_\nu(t-s)dW_s(\theta(r,\omega))
$$
holds for $a.e.\omega.$.
Notice at first, for $r>0$, let $s'=s+r$, then
\begin{eqnarray}
  u(t,\theta(r,\omega))
 &=&\frac{1}{2} \int_{-\infty}^{t+r} T_\nu(t+r-s')\frac{\partial
  u^2(s'-r,\theta(r,\omega))}{\partial x}ds' \nonumber\\
  &&
  \hskip1cm
  +\int_{-\infty}^{t+r}T_\nu(t+r-s')dW_{s'-r}
  (\theta(r,\omega))\nonumber\\
  &=&\frac{1}{2} \int_{-\infty}^{t} T_\nu(t+r-s')\frac{\partial
  u^2(s'-r,\theta(r,\omega))}{\partial x}ds'\nonumber\\
  &&
    \hskip1cm +\int_{-\infty}^{t}T_\nu(t+r-s')dW_{s'-r}
  (\theta(r,\omega))\nonumber\\
  & &+\frac{1}{2} \int_{t}^{t+r} T_\nu(t+r-s')\frac{\partial
  u^2(s'-r,\theta(r,\omega))}{\partial x}ds'\nonumber\\
  &&
   \hskip1cm  +\int_{t}^{t+r}T_\nu(t+r-s')dW_{s'-r}
  (\theta(r,\omega))\nonumber\\
  &=&T_\nu(r)\Big[\frac{1}{2}\int_{-\infty}^{t} T_\nu(t-s')\frac{\partial
  u^2(s'-r,\theta(r,\omega))}{\partial x}ds'\nonumber\\
  &&
  \hskip1cm   +\int_{-\infty}^{t}T_\nu(t-s')dW_{s'-r}
  (\theta(r,\omega))\nonumber\Big]\nonumber\\
  & &+\frac{1}{2} \int_{t}^{t+r} T_\nu(t+r-s')\frac{\partial
  u^2(s'-r,\theta(r,\omega))}{\partial x}ds'\nonumber\\
  &&
  \hskip1cm   +\int_{t}^{t+r}T_\nu(t+r-s')dW_{s'-r}
  (\theta(r,\omega)).\nonumber
\end{eqnarray}
On the other hand, let $y(\cdot,\omega)=u(\cdot-r,\theta(r,\omega))$, then we have
\begin{eqnarray}
  y(t+r,\omega)&=&T_\nu(r)\Big[\frac{1}{2}\int_{-\infty}^{t} T_\nu(t-s)\frac{\partial
  y^2(s,\omega)}{\partial x}ds+\int_{-\infty}^{t}T_\nu(t-s)dW_{s}
  (\omega)\nonumber\Big]\nonumber\\
  & &\ \ \ +\frac{1}{2} \int_{t}^{t+r} T_\nu(t+r-s)\frac{\partial
  y^2(s,\omega)}{\partial x}ds+\int_{t}^{t+r}T_\nu(t+r-s)dW_{s}
  (\omega)\nonumber\\
  &=&T_\nu(r)y(t,\omega)+\frac{1}{2} \int_{t}^{t+r} T_\nu(t+r-s)\frac{\partial
  y^2(s,\omega)}{\partial x}ds\nonumber\\
  &&
  \ \ \ +\int_{t}^{t+r}T_\nu(t+r-s)dW_{s}
  (\omega).\nonumber
\end{eqnarray}
By the uniqueness of equation (\ref{burgers u}) (see Theorem 2.6), we
know $y(t+r,\omega)=u(t+r,\omega)$ a.s. So,
$u(\cdot,\theta(r,\omega))=u(\cdot+r,\omega)$ a.s..

Using a similar argument, for $r<0$,
\begin{eqnarray}
  & &u(t,\theta(r,\omega))=\frac{1}{2} \int_{-\infty}^{t+r} T_\nu(t+r-s')\frac{\partial
  u^2(s'-r,\theta(r,\omega))}{\partial x}ds'  \nonumber \\
  & &\ \ \ \ \ \ \ \ \ \ \ \ \ \ \ \ \ \ \ \ +\int_{-\infty}^{t+r}T_\nu(t+r-s')dW_{s'-r}
  (\theta(r,\omega)).  \label{3.21}
\end{eqnarray}
Let $y(\cdot,\omega)=u(\cdot-r,\theta(r,\omega))$, then (\ref{3.21}) becomes
\begin{eqnarray*}
  y(t+r,\omega)&=&\frac{1}{2}\int_{-\infty}^{t+r} T_\nu(t+r-s')\frac{\partial
  y^2(s',\omega)}{\partial x}ds'+\int_{-\infty}^{t+r}T_\nu(t+r-s')dW_s.
\end{eqnarray*}
The uniqueness of (\ref{3.1}) proved in Lemma 3.1 shows that
$y(t+r,\omega)=u(t+r,\omega)$, a.s. i.e.
$u(\cdot,\theta(r,\omega))=u(\cdot+r,\omega)$ a.s.\ \hfill $\Box$
\medskip

Now, let's construct the solution of equation (\ref{3.1}). We denote for $n\in
{\bf Z}^+$
\begin{equation}u_n(t,x,\omega)=
  \left\{
     \begin{array}{ll}
        u(t,-n,0,\omega),& \ \ {\rm for \ }  t>-n,\\
        0,            & \ \ {\rm for } \ t\leq-n .
     \end{array}            \right.
  \label{burgers v}
\end{equation}
Then we can prove the following lemma.


\medskip

{\bf Lemma 3.3} {\it  For $a.e.\omega$, and any $N\in {\bf Z}^+$,
$u_n(t,x,\omega)\rightarrow u^*(t,x,\omega)$ in
$C([-N,N],L^2[0,1])$ as $n\rightarrow\infty$ under the norm
$|u(t)|_{\infty, L^2,N}\equiv\sup_{t\in [-N,N]}|u(t)|_2$. Moreover
$u^*$ satisfies
 (\ref{3.1}), (\ref{3.1c}) and (\ref{3.2}).}
\medskip

{\bf Proof}: At first, following Theorem A.7, we know that $u_n$ is
a Cauchy sequence in $C([-N,N], L^2[0,1])$ under the norm
$|u|_{\infty, L^2,N}$. Because
the space $C([-N,N], L^2)$ is complete, there is a $u^*$ such that
$\lim_{n\rightarrow\infty} u_n=u^*$ in $C([-N,N], L^2[0,1])$. Since
$N$ is arbitrary, $u^*(t,\omega)$ is defined for all time.

Secondly, for any $t$ and $t_0<t$, we will show that $u^*$ satisfies,
\begin{equation}
  u^*(t)=T_\nu(t-t_0)u^*(t_0)+\frac{1}{2}\int_{t_0}^{t} T_\nu(t-s)\frac{\partial
  (u^*(s,\omega))^2}{\partial x}ds+\int_{t_0}^{t}T_\nu(t-s)dW_{s}.
  \label{3.30}
\end{equation}
For this, similar to the proof of formula (4.13) in \cite{[Fla]}, for any $N$, any sufficiently
large $n$ and $a.e.\omega$, we have
\begin{eqnarray*}
 \int_{-N}^N \Big|\frac{\partial u_n(s)}{\partial
 x}\Big|^2_2ds<\xi_N(\omega)<\infty.
\end{eqnarray*}
This means that we can find a subsequence, still denoted by $u_n$,
weakly converge to $u^*$ in $L([-N,N],H^1_0)$, by Alaoglu
Compactness Theorem. Therefore for $a.e.\omega$, for any $N$,
$u^*\in L([-N,N],H^1_0)$. Moreover, using the same estimate as in
Lemma 2.7, we get
\begin{eqnarray*}
  \frac{1}{2}\int_{t_0}^{t} T_\nu(t-s)\frac{\partial
  (u^*(s,\omega))^2}{\partial x}ds\in L^2[0,1].
\end{eqnarray*}
Thus
\begin{eqnarray*}
   & &\Big|\int_{t_0}^{t} T_\nu(t-s)\frac{\partial
  (u^*(s,\omega))^2}{\partial x}ds-\int_{t_0}^{t} T_\nu(t-s)\frac{\partial
  (u_n(s,\omega))^2}{\partial x}ds\Big|_2^2\nonumber\\
   &=&\int_0^1\Big[\int_{t_0}^t\int_0^1p(t-s,y,x)\Big(\frac
      {\partial (u^*(s,\omega))^2}{\partial x}-\frac{\partial (u_n(s,\omega))^2}{\partial
      x}\Big)dxds\Big]^2dy\nonumber\\
      &\leq &C\sup_{t_0\leq s\leq t}[|u_n(s)+u^*(s)|_2^2]\int_{t_0}^t
 \frac{1}{(t-s)^{\frac{3}{4}}}|u_n(s)-u^*(s)|_2^2ds,
 \label{3.31}
\end{eqnarray*}
using the same estimate technique in the proof of the uniqueness of Theorem
2.6. Since $u_n\rightarrow u^*$ in $C([-N,N], L^2[0,1])$ under norm
$|u(t)|_{\infty, L^2,N}$, it is easy to know that
\begin{eqnarray*}
   & &\Big|\int_{t_0}^{t} T_\nu(t-s)\frac{\partial
  (u^*(s,\omega))^2}{\partial x}ds-\int_{t_0}^{t} T_\nu(t-s)\frac{\partial
  (u_n(s,\omega))^2}{\partial x}ds\Big|_2^2\rightarrow
0.
\end{eqnarray*}
 At the same time, obviously $u_n(t)$, $T_\nu(t-t_0)u_n(t_0)$
strongly converge to $u^*(t)$ and $T_\nu(t-t_0)u^*(t_0)$ in
$L^2[0,1] $ respectively, hence (\ref{3.30}) holds.

Due to Theorem A.4, $\sup_n\sup_{t\in {\bf R}}E|u_n(t)|^{2p}_2<\infty$, for any $
p\geq 1$. This implies that
\begin{equation}
  \sup_{t\in {\bf R}}E|u^*(t)|^{2p}_2<\infty,\ p\geq 1. \label{3.33}
\end{equation}

Finally, let's prove that $u^*$ satisfies (\ref{3.1}). From (\ref{3.30}), it is
easy to know, for any $0<m<n$,
\begin{eqnarray*}
  &&\frac{1}{2}\int_{-n}^{-m} T_\nu(-s)\frac{\partial
  (u^*(s,\omega))^2}{\partial
  x}ds\\
  &=&-\int_{-n}^{-m}T_\nu(-s)dW_s-T_\nu(n)u^*(-n)+T_\nu(m)u^*(-m).
\end{eqnarray*}
Thus
\begin{eqnarray*}
  E\Big|\int_{-n}^{-m} T_\nu(-s)\frac{\partial
  (u^*(s,\omega))^2}{\partial
  x}ds\Big|^2_2&\leq&C\Big[E\Big|\int_{-n}^{-m}T_\nu(-s)dW_s\Big|^2_2\\
  &&\ \
+  E\Big|T_\nu(n)u^*(-n)\Big|_2^2+E\Big|T_\nu(-m)u^*(m)\Big|_2^2\Big]\\
  &=&I+II+III.
\end{eqnarray*}
It is easy to know that $I\rightarrow 0$ as $n,m\rightarrow\infty$. Moreover,
by Poincar\'{e}'s inequality (see \cite{chen}), we know $II\leq e^{-\nu\lambda_1n}E|u^*(-n)|_2^2$
and $III\leq e^{-\nu\lambda_1m}E|u^*(-m)|_2^2$, then from (\ref{3.33}) we know that $II$
and $III$ converge to $0$ as $n,m\rightarrow\infty$. So, there is a subsequence
$n'$ such that, for $a.e.\omega$, $\int_{-n'}^0T_\nu(-s)\frac{\partial
  (u^*(s,\omega))^2}{\partial
  x}ds\rightarrow \int_{-\infty}^0T_\nu(-s)\frac{\partial
  (u^*(s,\omega))^2}{\partial
  x}ds$ and
$\int_{-n'}^0T_\nu(-s)dW_s\rightarrow
\int_{-\infty}^0T_\nu(-s)dW_s$ in $L^2[0,1]$ as $n'\rightarrow\infty$,  and
$T(n)u^*(-n)\rightarrow 0$ in $L^2[0,1]$ as $n\rightarrow\infty$.
Therefore, $u^*$ satisfies equations (\ref{3.1}) and (\ref{3.1c})
and Lemma 3.2 implies (\ref{3.2}) holds.\ \ \hfill $\Box$

\medskip

{\bf Remark}:\ Assume that $\tilde{u}$ is another stationary
solution for (\ref{burgers u}), since $\theta(t)$ is a
$P$-preserving ergodic Wiener shift on $\Omega$, $\tilde{u}(t)$
have the same law for any $t$, and its law is the invariant
measure for (\ref{burgers u}). Similar to the proof of ergodicity
of stochastic Burgers equation in Chapter 14 of \cite{[DPZ2]}, we can
conclude that the invariant measure of (\ref{burgers u}) is
unique. Therefore, the law of $\tilde{u}(t)$ identifies that of $u^*(t)$.
So, $\tilde{u}$ satisfies the moment estimate (\ref{3.1c}).
Because $\tilde{u}$ is the solution for  (\ref{burgers u}) for any
$t\in{\bf R}$, using the same reasoning as inequalities of
(\ref{3.1d}), (\ref{3.1e}) in the proof Lemma 3.1, we know that
$u^*=\tilde u$ so there
exists unique stationary solution for  (\ref{burgers u}).

\vskip5pt

Now we have proved that equation (\ref{3.1}) has a unique solution $u^*$
in $C((-\infty,+\infty),L^2([0,1])$ and for any $r$, (\ref{3.2}) holds almost surely. On the
other hand, the stationary solution of (\ref{burgers u}) is unique so must satisfy (\ref{3.1})
and ({\ref{3.1c}). It remains to prove (\ref{3.2}) holds for all $r$ almost surely by a perfection
argument.
\vskip5pt


{\bf Theorem 3.4}{\it \ \ Under the Condition A, there exists
a unique stationary solution $\hat{u}^*(\cdot)$ for the stochastic
Burgers equation (\ref{burgers u}) satisfying:-
\begin{enumerate}
  \item[(i)] $\hat{u}^*$ is $({\cal B}({\bf R})\otimes{\cal F},{\cal
  B}(L^ 2[0,1]))$ measurable;
  \item[(ii)] $\hat{u}^*(\cdot,\omega)$:${\bf R}\rightarrow L^2[0,1]$
  is continuous for all $\omega \in \Omega$;
  \item[(iii)] For all $r,t\in {\bf R}$ and $\omega \in \Omega$,
  $\hat{u}^*(t, \theta(r,\omega))=\hat{u}^*(t+r,\omega)$;
  \item[(iv)] Let $\tilde{N}\equiv$ $\{\omega:\hat{u}^*(t,\omega)\neq u^*(t,\omega)\  \mbox{\it for some}\  t\in {\bf
  R}\}$, then $\tilde{N}\in {\cal F}$ and $P(\tilde{N})=0$. This
  implies that $\hat{u}^*$ and $u^*$ are indistinguishable and satisfies (\ref{3.1}) and (\ref{3.1c});
  \item[(v)] $\hat{u}^*$ is $({\cal F}_t, \ t\in {\bf R})$
  adapted.
\end{enumerate}

Moreover, let $Y(\omega)=\hat{u}^*(0,\omega)$, which is ${\cal
F}_0$ measurable. Then for all $\omega\in \Omega$,
\begin{equation}
  Y(\omega)=\frac{1}{2}\int_{-\infty}^0T_\nu(-s)\frac{\partial
  (Y(\theta(s,\omega))^2)}{\partial
  x}ds+\int_{-\infty}^0T_\nu(-s)dW_s(\omega),\label{3.4}
\end{equation}
thus, $Y(\omega)$ is the stationary point and
$\hat{u}^*(t,0,Y(\omega),\omega)=Y(\theta(t,\omega))$ for $t\geq
0$.}
\vskip5pt

{\bf Proof}:  The above remark implies the uniqueness.
Now, we use a similar argument of the perfection of crude
cocycle in \cite{[Ar]} (see p15) or \cite{[KS]} to prove (i)--(v) and (\ref{3.4})
hold. First note due to Lemmas 3.2 and 3.3, we have for any
$r\in {\bf R}$, there is a $N_r\in {\cal F}$ such that $P(N_r)=0$ and for any $\omega\in N_r^c$,
$$
 u^*(t, \theta(r,\omega))=u^*(t+r,\omega),\ \mbox{\it for any }\ t.
$$
Denote
\begin{eqnarray*}
 M&\equiv&\{(r,\omega)\in {\bf R}\times\Omega,\ u^*(t,
 \theta(r,\omega))=u^*(t+r,\omega)\ \mbox{\it for any }\ t\};\\
 \Omega_0&\equiv&\{\omega\in \Omega,\ (s,\omega)\in M,\ \lambda-a.e.\
 \ s\in {\bf R}\};\\
 \Omega_1&\equiv&\{\omega\in \Omega,\ \theta(s,\omega)\in \Omega_0,\
\lambda-a.e.\ \ s\in {\bf R}\}.
\end{eqnarray*}
Here $\lambda(dx)$ is the Lebesgue measure on ${\bf R}$. We will prove the theorem in the following 7 steps.

Step 1. We should show that $M$ is a measurable set in
${\cal B}({\bf R})\otimes{\cal F}$.

Since $u^*(t,\omega)$ is continuous in $C((-\infty,\infty), L^2)$,
we have
$$
  M=\cap_{t\in {\bf Q}}\{(r,\omega)\in {\bf R}\times\Omega,\ u^*(t,
 \theta(r,\omega))=u^*(t+r,\omega)\},
$$
where $\bf Q$ is the set of rational number in $\bf R$. By the
construction of $u^*$, we know that $u^*(t,\omega)$ is $({\cal
B}({\bf R})\otimes{\cal F},{\cal
  B}(L^ 2[0,1]))$ measurable. Because $\theta(r,\omega)$ is ${\cal
B}({\bf R})\otimes {\cal F}$ measurable, these imply that
$u^*(t+r,\omega)$ and $u^*(t,\theta(r,\omega))$ both are $({\cal
B}({\bf R})\otimes{\cal F},{\cal
  B}(L^ 2[0,1]))$  measurable. Furthermore, since $L^2[0,1]$
is an Hausdorff and second countable space, we have  that $
\{(r,\omega)\in {\bf R}\times\Omega,\ u^*(t,
 \theta(r,\omega))=u^*(t+r,\omega) \}
$ is $({\cal B}({\bf R})\otimes{\cal F},{\cal
  B}(L^ 2[0,1]))$ measurable for every $t \in {\bf Q}$, therefore,
  $M$ is $({\cal B}({\bf R})\otimes{\cal F},{\cal
  B}(L^ 2[0,1]))$ measurable.

Moreover, using Fubini's theorem, we get $\lambda\otimes P({\bf
R}\times\Omega\setminus M)=0$.

{Step 2}. \ Taking
$\nu(dt)\equiv\frac{1}{\sqrt{2\pi}}e^{-\frac{t^2}{2}}\lambda(dt)$,
it is easy to know
that $
  \omega\in \Omega_0$ iff $\int_{\bf R}{\bf
  1}_M(s,\omega)\nu (ds)=1.
$
Fubini's theorem implies $\Omega_0\in {\cal F}$, and moreover $\lambda \otimes (R\times \Omega\setminus M)=0$, one can see that
$$
 \int_{\bf R}\int_\Omega{\bf
  1}_M(s,\omega)P(d\omega)\nu(ds)=1\ \ {i.e}\ \ \int_\Omega\int_{\bf R}{\bf
  1}_M(s,\omega)\nu(ds)P(d\omega)=1,
$$
so, $\int_{\bf R}{\bf
  1}_M(s,\omega)\nu (ds)=1\ P. a.s.$. This shows that $P(\Omega_0)=1$.
  Similarly, one can show $\Omega_1\in {\cal F}$ and
  $P(\Omega_1)=1$.

Now, we show that $\Omega_1$ is invariant under $\theta$, {\it
i.e.}, $\omega\in \Omega_1$ implies $\theta(t,\omega)\in \Omega_1$
for all $t\in {\bf R}$. For this, for $\omega\in \Omega_1$, there exists a
$\Lambda_\omega\in {\cal B}({\bf R})$ such that
$\lambda(\Lambda_\omega)=0$ and for any $ s\in \Lambda_\omega^c$,
$\theta(s,\omega)\in \Omega_0$. For any $t\in {\bf R}$, noticing
$\Lambda_{\theta(t,\omega)}=\Lambda_\omega-t$, it is obvious
$\lambda(\Lambda_{\theta(t,\omega)})=0$ and
$$
 \theta(s,\omega)=\theta(s-t+t,\omega)=\theta(s-t,\theta(t,\omega))\in
 \Omega_0.
$$
This implies that $\theta(t,\omega)\in \Omega_1$.

{Step 3}.\  Define $\hat{u}^*$,
\begin{equation}
 \hat{u}^*(t,\omega)\equiv
  \left\{
    \begin{array}{ll}u^*(t-s, \theta(s,\omega)), & {\rm if}\  \omega\in
      \Omega_1, \ \theta(s,\omega)\in \Omega_0\\
      {\bf 0}, &{\rm if} \ \omega\in \Omega_1
    \end{array}          \right.
\end{equation}
where $\bf 0$ is the constant function $0$ in $L^2[0,1]$.

Firstly, we should show $\hat{u}^*$ is well-defined. Assume that
$\omega\in \Omega_1, \theta(s,\omega)\in \Omega_0$ and
$\theta(u,\omega)\in \Omega_0$. Note there
exist $\Lambda'_{\theta(s,\omega)}$ and
$\Lambda'_{\theta(u,\omega)}\in {\cal B}({\bf R})$ such that
$\lambda(\Lambda'_{\theta(s,\omega)})=\lambda(\Lambda'_{\theta(u,\omega)})=0$,
and for any $\alpha\in \Lambda^{'c}_{\theta(s,\omega)}, \beta\in
\Lambda^{'c}_{\theta(u,\omega)}$, and for all $t\in R$,
$$
 u^*(t+\alpha,\theta(s,\omega))=u^*(t,\theta(\alpha,\theta(s,\omega)))=u^*(t,
 \theta(\alpha+s,\omega)),
$$
$$
 u^*(t+\beta,\theta(u,\omega))=u^*(t,\theta(\beta,\theta(u,\omega)))=u^*(t,
 \theta(\beta+u,\omega)).
$$
We can always take $\alpha\in \Lambda^{'c}_{\theta(s,\omega)},
\beta\in \Lambda^{'c}_{\theta(u,\omega)}$ such that
$\alpha+s=\beta+u$, so we have, for any $t$
\begin{eqnarray*}
 u^*(t+\alpha,\theta(s,\omega))=u^*(t+\beta,\theta(u,\omega))&=&u^*(t+\alpha+s-s,\theta(s,\omega))\\
 &=&u^*(t+\beta+u-u,\theta(u,\omega)).
\end{eqnarray*}
This implies that
$u^*(t-s,\theta(s,\omega))=u^*(t-u,\theta(u,\omega))$ for any
$t\in {\bf R}$.

{Step 4}.\  Obviously, $\hat{u}^*(\cdot,\omega): {\bf
R}\rightarrow L^2[0,1]$ is continuous for all $\omega\in \Omega$.

{Step 5}.\ We should show that $\hat{u}^*$ satisfies
for all $ r,t\in {\bf R}$, $\hat{u}^*(t,
\theta(r,\omega))=\hat{u}^*(t+r,\omega)$, for all $\omega \in
\Omega$.
The assertion is clear for $\omega\notin\Omega_1$.  We assume
$\omega\in \Omega_1$, thus $\theta(r,\omega)\in \Omega_1$ for any
$r\in {\bf R}$.
By the definition of $\hat{u}^*$, we have
$$
 \hat{u}^*(t+r,\omega)=u^*(t+r-s,\theta(s,\omega))
 \ \ {\rm for \ all }\ s\in R\setminus \Lambda^{\prime}_{\omega},\ t,r\in R,
$$
and
\begin{eqnarray*}
 \hat{u}^*(t,\theta(r,\omega))&=&u^*(t-u,\theta(u,\theta(r,\omega))=u^*(t-u,\theta(u+r,\omega))\\
  &&
  {\rm for \ all }\ u\in R\setminus \Lambda^{\prime}_{\theta (r,\omega)},\ t,r\in R.
 \end{eqnarray*}
On the other hand, note for any $r$, there exists $s\in R\setminus
\Lambda^{\prime}_{\omega}$, $u\in R\setminus
\Lambda^{\prime}_{\theta (r,\omega)}$ satisfying $r+u=s$. For if
this is not true, then there exists a $r\in R$ such that for any $
s\in R\setminus \Lambda^{\prime}_{\omega}$, $u\in R\setminus
\Lambda^{\prime}_{\theta (r,\omega)}$, $r+u=s$ cannnot be true.
That is to say that for any $u\in R\setminus
\Lambda^{\prime}_{\theta (r,\omega)}$, $r+u\notin R\setminus
\Lambda^{\prime}_{\omega}$ so $r+u\in \Lambda^{\prime}_{\omega}$.
This is certainly not true since $\nu (R\setminus
\Lambda^{\prime}_{\theta (r,\omega)})=1$ and $\nu
(\Lambda^{\prime}_{\omega})=0$. Thus the assertion holds for any
$t$ and $r$.

{Step 6}.\ Let
\begin{equation}
 B(s,t,\omega)\equiv
  \left\{
    \begin{array}{ll}u^*(t-s, \theta(s,\omega)) & \ \ \omega\in
      \Omega_1, \ \theta(s,\omega)\in \Omega_0\\
      {\bf 0}&  \ \ \omega\in \Omega_1
    \end{array}          \right.
\end{equation}
Using the same reasoning as the step 6 on p20 of \cite{[Ar]}, we know that
$\hat{u}^*(t,\omega)$ is ${\cal B}({\bf R})\otimes{\cal F}$
measurable.

Because $\hat{u}^*$ and $u^*$ both are continuous and ${\cal
B}({\bf R})\otimes{\cal F}$ measurable,
$$
 \{\omega: \hat{u}^*(t,\omega)=u^*(t,\omega),\ \ {\rm for \ all}\  t\in {\bf
 R}\}=\cap_{t\in {\bf Q}}\{\omega:
 \hat{u}^*(t,\omega)=u^*(t,\omega)\}\in {\cal F}.
$$
Moreover, for any $t$, and $ \omega\in \Omega_0\cap\Omega_1$,
we obtain $\hat{u}^*(t,\omega)=u^*(t,\omega)$. Since
$P(\Omega_0\cap\Omega_1)=1$, we know
$P(\{\omega:\hat{u}^*(t,\omega)=u^*(t,\omega)\})=1$. All of these
imply that $\hat{u}^*$ and $u^*$ are indistinguishable.

{Step 7}.\ We should prove $\hat{u}^*$ is $({\cal F}_t)_ {t\in
{\bf R}}$ adapted. Due to the construction of $u^*$, we know that
$u^*$ is adapted {\it i.e.} $u^*(t,\cdot)\in {\cal F}_t$. (see the
beginning of section 2 for the definition of ${\cal F}_t$) It is
easy to know that $B(s,t,\omega)\in {\cal F}_t$ for any $s,t$. By
(\ref{theta}) and $P(\Omega_0)=1$, $P(\Omega_1)=1$, this means
that $\hat{u}^*$ is adapted to $({\cal F}_t)_ {t\in {\bf R}}$. \ \
\hfill $\Box$

 \vskip5pt

Furthermore, due to Theorem 2.10, Theorem A.5 in this paper and Theorem 2.1.1, 2.1.2
and 2.2.1 in \cite{[MZZ]},
we have the following dynamical behaviour near $Y(\omega)$. \vskip5pt

{\bf Theorem 3.5} \ {\it If condition A holds, then
 for solution of (\ref{burgers u}), $u(t,0,u_0,\omega)$ stating at time $0$,  with initial data $u_0\in L^2[0,1]$, we have,
\begin{equation}
  lim_{t\rightarrow \infty}{1\over t}\log|u(t,0,u_0,\omega)-Y(\theta(t,\omega))|_2\leq
  -\delta, \ \ P.a.s.
\end{equation}
and
\begin{equation}
 lim_{t\rightarrow \infty}\log\|Du(t,Y(\omega),\omega)\|_{L(L^2[0,1]))}\leq
 -\delta,\ \ \forall \omega\in \Omega,
\end{equation}
where $L(L^2[0,1])$ is the Hilbert space of all bounded linear
operators with operator norm, and $\delta \in (0,\delta _0)$. }
\vskip5pt

The theorem says that $(u,\theta)$ is a $C^1$ perfect cocycle on
$L^2[0,1]$, and the stationary point $Y(\omega)$ is hyperbolic,
especially, its the largest Lyaponov spectrum not larger than $-\delta(<0)$ and $L^2[0,1]$ is the stable manifold of $Y(\omega)$
for all $\omega\in \Omega^*.$


\section*{Appendix A}

Although we believe that experts in this field are familiar
with some properties here and some of them can be found in \cite{[DPZ1]},
we would like to include them here for completeness.

It is easy to see from Theorem 5.4 in \cite{[DPZ1]} that $\w$ is the unique weak
solution in $L^2[0,1]$ of the following SPDE,
\begin{equation}
  \left\{
     \begin{array}{ll}
       dX(x,t)=\nu\Delta X(x,t)dt+dW(t), \ t\geq 0; x\in (0,1),\\
       X(0,t)=X(1,t)=0,\\
       X(x,0)=0, \ \ x\in (0,1).\ ,
     \end{array}           \label{1} \right.
  \label{burgers}
\end{equation}
Therefore for any $e_k$, by integration by parts
$$\langle\w(t),e_k\rangle=\int_0^t\langle \w(s), \nu\Delta e_k\rangle ds+\sigma_k(B_k(t)-B_k(0)).$$
Let $\wk=\langle\w(t),e_k\rangle$. Obviously, $\wk$ satisfies the
following O-U equation in 1 dimension:
\begin{equation}
 \wk(t)=-\nu\pi^2k^2\int_0^t\wk(s)ds+\sigma_k(B_k(t)-B_k(0)),
\end{equation}
so
\begin{equation}
  \wk(t)=\sigma_k\int_0^t e^{-\nu\pi^2k^2(t-s)}dB_k(s).\label{wk}
\end{equation}
Let $\w^{(n)}(t)\equiv P_n\w(t)=\sum_{k=1}^n\wk(t)e_k$, then by It\^{o}'s
formula, we have
\begin{equation}
|\w^{(n)}(t)|^2_2+\nu\int_0^t\Big|\p\w^{(n)}(s)\Big|_2^2ds=\sum_{k=1}^n\int_0^t
\wk(s)\sigma_kdB_k(s)+t\sum_{k=1}^n\sigma_k^2.\label{I}
\end{equation}

\vskip5pt

{\bf Lemma A.1:} {\it $\displaystyle \sum_{k=1}^n\int_0^t
\wk(s)\sigma_kdB_k(s)$ is a martingale.}

\proof In fact, we only need to show that
\begin{equation}
  E\int_0^t \Big(\wk(s)\sigma_k\Big)^2ds<\infty .\label{E}
\end{equation}
In fact from  (\ref{wk}), we have
\begin{eqnarray*}
E\int_0^t \Big(\wk(s)\sigma_k\Big)^2ds&=&\sigma _k^4E\int_0^t\Big(\int_0^se^{-\nu\pi^2k^2(s-r)}dB_k(r)\Big)^2ds\\
 &=&\sigma _k^4\int_0^t E\Big(\int_0^se^{-\nu\pi^2k^2(s-r)}dB_k(r)\Big)^2ds\\
 &=&\sigma _k^4\int_0^t\int_0^se^{-2\nu\pi^2k^2(s-r)}drds\\
 &=&\frac{\sigma _k^4}{2\nu\pi^2k^2}\Big(t+\frac{1}{2\nu\pi^2k^2}(e^{-2\nu\pi^2k^2t}-1)\Big)<\infty
\end{eqnarray*}
\endpf


\vskip5pt

Using Lemma A.1, (\ref{I}) and $\sum_{k=1}^\infty
\sigma_k^2<\infty$, one can get it immediately \vskip5pt


{\bf Theorem A.2:} {\it $\displaystyle
E\Big(\int_0^t\Big|\p\w(s)\Big|^2_2ds\Big)<\infty$, therefore, for any $t$,
$\w\in L^2([0,t], H^1_0)\ \ a.s.$.}

\vskip5pt


By theorem 16.2 in \cite{[LSU]}  (p413), one easily deduces that
\vskip5pt

{\bf Proposition A.3}: \ {\it There exist positive constants $c_1, c_2$ such that
$$
  \Big|\frac{\partial p(t,x,y)}{\partial y}\Big|\leq
  \frac{c_1}{t}e^{-\frac{(x-y)^2}{2c_2t}},
$$
for all $t>0,\ x,y\in [0,1]$. Then one can pick a positive constant $c_3$ such that
$$
 \int_0^1\frac{c_3}{t}e^{-\frac{y^2}{2c_2t}}dy\leq 1.
$$
}


The following results are from \cite{[Ma]}. Although he
dealt with the case of 2D stochastic Navier-Stokes equation, all of his
techniques and reasonings can be used to 1D stochastic Burgers
equations with Dirichlet boundary condition. So, we omit the proof of the
theorems.
\vskip5pt

{\bf Theorem A.4} (see Lemma 2 in \cite{[Ma]} for details)\ {\it Assume the
initial condition satisfies $E|u(t_0)|_2^{2p}<\infty$ for some $p\geq 1$, then
there is a constant $c$ such that
$$E|u(t,t_0, u(t_0))|_2^{2p}\leq
c^p(p-1)!.$$ }
\vskip5pt

{\bf Theorem A.5} (see Theorem 1 in \cite{[Ma]})\ {\it Assume Condition A holds.
Fix a $\delta\in (0,\delta_0)$ and a time $t_0$, Let $u_0\in L^2$ be an initial
condition, measurable with respect to ${\cal F}_0$, such that
$E|u_0|_2^{2p}<\infty$ for some $p>1$. Let
$\overline{u}\equiv\overline{u}(t,t_0,\overline{u}_0)$ denote  the
solution starting from some other arbitrary initial data $\overline{u}_0\in
L^2$. Then, there is a positive integer-valued random time
$\tau(\delta,t_0,u_0)$, independent of $\overline{u}$, such that
$$
 |u(t)-\overline{u}(t)|_2^2\leq|u_0-\overline{u}_0|^2_2 e^{-2\delta(t-t_0)}
$$
for all $t>t_0+\tau$. In addition, $E(\tau^q)<\infty$ for any $q\in (0,p-1)$.
}

\vskip5pt

{\bf Theorem A.6} (see Theorem 2 in \cite{[Ma]})\ {\it Assume Condition A holds.
Fix a $\delta\in (0,\delta_0)$ and a $t\in {\bf R}$. Let $\{u_0(n)\}$ be a
sequence of random variable with $n\in {\bf Z}^+$. Assume that the $u_0(n)$ are
measurable with respect to ${\cal F}_{t-n}$ and that $E|u_0(n)|_2^{2p}$ is
uniformly bounded in $n$ for some $p>2$. Then the following hold:

(1)\  With probability one,there exists a random ${\bf Z}$-valued time
$\overleftarrow{n}(\epsilon,\delta,t,\omega)>0$ such that for real $s>0$ and
all $n\in{\bf Z}$ with $n>\overleftarrow{n}$ we have
$$
  \sup_{u'_0 \in
  A_n}|u(t+s,t-n,u_0(n),\omega)-u(t+s,t-n,u'_0,\omega)|_2\leq\epsilon\delta^2|n|e^{-\delta(|n|+s)}.
$$
Here $A_n$ is the set $\{u'_0: |u'_0|^2_2\leq \frac{\epsilon\delta^2}{2}|n|\}$.
In addition, $E(\overleftarrow{n}^q)<\infty$ for any $q\in(0, p-2)$.

(2)\ Let $\{\overline{u}_0(n)\}$ be a second sequence of random variables with
$n\in {\bf Z}^+$ measurable with respect to ${\cal F}_{t-n}$ and
$E|\overline u_0(n)|_2^{2p}$ uniformly bounded in $n$ for some $p>2$. Then with probability one,
there exists
another ${\bf Z}$-valued random time $\overleftarrow{n'}$ such that for real $s>0$ and all $n\in {\bf Z}$ with
$n>\overleftarrow{n'}$ we have
$$
 |u(t+s,t-n,u_0(n))-u(t+s,t-n,\overline{u}_0(n))|_2\leq
 \epsilon\delta^2|n|e^{-\delta(|n|+s)}.
$$
Again, $E((\overleftarrow{n'})^q)<\infty$ for any $q\in(0,p-2)$.}
\vskip5pt

{\bf Theorem A.7} (see Corollary 1 in \cite{[Ma]})\ {\it Under Condition A, fix
$t\in {\bf Z}$ and a $\delta\in (0,\delta_0)$. Given any $\epsilon>0$, with probability one,
there is
a positive ${\bf Z}$-valued random time $n^*(\epsilon,\delta,t_1)$ such that
 for all $\tau\leq 0$ and all $n_1, n_2\in {\bf Z}$, if
$n_1,n_2<t-1-n^*$, we have
$$
 |u(t+\tau,n_1,0)-u(t+\tau,n_2,0)|_2\leq \epsilon e^{-\delta\tau}.
$$
Furthermore, $n^*(\omega)$ is a stationary random variable with
all moments finite.}

\medskip
{
{\bf  Acknowledgements}\ \ One of the authors, YL,
would like to acknowledge the financial support of
EPSRC GR/R69518 and partial supports of NSFC
(No.10531070, 10171101, 10101002), SRF for ROCS. Both
authors would like to thank Prof. D. Elworthy for his invitation to
visit the University of Warwick, and Prof. S. Peng for his invitation to visit
Shangdong University.

\end{document}